\documentclass[twoside,10pt]{article}
\usepackage{ITMM2024}
\usepackage{caption}
\usepackage{subcaption}
\newcommand{\Ii}{ \mathbf{1} }
\begin{document}
\English
\title[Information Spreading in 
Random Graphs Evoving by Norros-Reittu 
Model]
{Information Spreading in Random Graphs Evolving by Norros-Reittu 
Model}{Information Spreading in Random Graphs Evoving by Norros-Reittu 
Model}
\author[N.\, M.~Markovich, D.\,V.~Osipov]
       {N.\,M.~Markovich\affilnote{1}, D.\,V.~Osipov\affilnote{2}}
       {Markovich~N.\, M., Osipov~D.\, V.}
\organization{\affilnote{1}V.A. Trapeznikov Institute of Control Sciences\\
Russian Academy of Sciences \\
Profsoyuznaya Str. 65, 117997 Moscow 
Russia 
\\
\affilnote{2}
The Department of Discrete Mathematics \\
Moscow Institute of Physics and Technology, Moscow
Russia }%

\maketitle

\grant{Markovich N.M. was supported by the Russian Science Foundation RSF, project number 24-21-00183.}

\begin{abstract}
The paper is devoted to the spreading of a message within the random graph evolving by the Norros-Reittu preferential attachment model. The latter model forms random Poissonian numbers of edges between newly added nodes and existing ones. For a pre-fixed time $T^*$, the probability mass functions of the number of nodes obtained the message and the total number of nodes in the graph, as well as the distribution function of their ratio are derived. To this end, the success probability to disseminate the message from the node with the message to the node without message is proved. The exposition is illustrated by the simulation study.

{\bf Keywords:} 
Evolution, random graph, Norros-Reittu preferential attachment model, information spreading
\end{abstract}

\section*{Introduction}
Preferential attachment models (PAMs) attract the attention of many researchers due to numerous applications to the evolution of real-world networks, such as social ones. PAMs allow us to model heavy-tailed empirical distributions of in- and out-degrees of nodes, as well as their PageRanks. An evolution starts from an initial graph $G_0$ that contains at least one isolated node. Usually, it is assumed that a new node is attached to existing nodes 
with the prefixed $m\ge 1$ new edges. However, PAMs became more flexible and realistic when $m$ is a random number. In the Norros-Reittu model 
\cite{Norros} and Poisson PAM in \cite{Wang} the number of new edges is modeled by the Poisson distribution. The randomness of $m$ is especially important during the start-up of the network. The discrepancy in the degree distribution between the traditional and Poisson PAMs becomes negligible over time, \cite{Wang}.
\par
Spreading  information is an important problem in random networks such as multi-agent systems, telecommunication overlay networks, parallel computation in  computational Grid systems \cite{Mosk-Aoyama2006}, \cite{Censor-Hiller2010}, social networks and spread of infections \cite{Pat}, \cite{Lang}, \cite{Pang} and gossip algorithms \cite{Shah}.
\par 
Let   $G_n=\left( V_n, E_n\right)$, $n=0,1,2, \dots$  denote the sequence of random graphs evolved by an evolution model, where $V_n$ is a set of nodes and $E_n$ is a set of edges. We denote the cardinality of the set $A$ as $\# A$. The evolution begins from the initial graph $G_0$. We may assume that $G_0$ contains a single node with the message.
\par
Our objective is to study a spreading of a unique message in   undirected graphs evolved by the Norros-Reittu model. We aim to analyze the growth of the set of nodes $S_k$ that obtained the message at the evolution step $k$.  

Let $N_k = \#V_k$ be the number of nodes in the graph in the evolution step $k$. By convention, we assume that the message can be transmitted from a node $w \in S_k$ to a node $v \in V_{k+1} \setminus S_k $ at step $k+1$ if there exists an edge between these nodes. The transmission of the message is governed by the ticks of a Poissonian clock, which drive the evolution of the graph, followed by the propagation of the message to the new node. We assume that message transmissions occur for a fixed time $T^*$. $K^*$ is the number of ticks of the clock (or evolution steps) up to $T^*$.

We aim to consider the dynamics of $N_k$ and $S_k$ over time and obtain distributions of $\#S_k$, $N_k$, $\#S_{K^*}$, $N_{K^*}$ , $\#S_{K^*} / N_{K^*}$.

\par
The paper is organized as follows. The statement of the problem is formulated in Section \ref{Sec1}. In Section \ref{Sec2} related results are recalled.  Our main results are presented in Section \ref{Sec3}. The simulation study is given in Section \ref{Sec4}. We finalize with conclusions.

\section{Statement of the problem}\label{Sec1}
We aim to obtain the following probability mass functions (pmfs):
\begin{eqnarray}\label{2}&& P\{\#S_{K^*}=i| \Lambda\}, \qquad P\{N_{K^*}=i\},
\end{eqnarray}
and  the distribution function of the proportion of nodes with the message at maximum step $K^*$ (a tick of a Poissonian clock) in the fixed time $T^*$
\begin{eqnarray}\label{3a}&& P\{\#S_{K^*}/N_{K^*}\le x| \Lambda\}.\end{eqnarray}
A sequence $\Lambda$ is defined in the next section.
A new node appends to  the graph by a tick of a Poissonian clock. The interarrival times between  consecutive pairs of new nodes are exponentially distributed with intensity $\lambda$. The probability that the number of  ticks $\nu(t)$ in time $t$ is equal to $k=0,1,2,...$ is
\begin{eqnarray}\label{3}
P\{\nu(t)=k\}=\frac{(\lambda t)^ke^{-\lambda t}}{k!}.
\end{eqnarray}
The number of evolution steps in time $T^*$ is $K^*=\nu(T^*)$, and its expectation is $E(K^*)=\lambda T^*$.

\section{Related work}\label{Sec2}
Here we consider the Norros-Reittu model \cite{Norros}, which is one of the key PAMs. In this model, new nodes join existing ones with probabilities depending on their mean degrees $\{\Lambda_i\}$, $i=0,1,...$ called $"$capacities$"$ by a random number of edges and create self-loops. This allows us to model more complex network structures due to possible multiple edges and self loops.

In \cite{Norros}, a sequence 
$\Lambda = (\Lambda_0, \Lambda_1, \ldots)$ is formed beforehand. Here,  $\{\Lambda_i\}$ are assumed to be independent, strictly positive  random variables (r.v.s) with a common distribution function $F$ such that
\begin{eqnarray}\label{0}P\{\Lambda_1>x\}&= & x^{-\tau +1}\ell(x), \qquad \tau>1,
\end{eqnarray}
$\ell(x)$ is a slowly varying function, i.e. by definition $\lim_{x\to\infty}\ell(tx)/\ell(x)=1$ for any $t>0$.
 For $\tau\in(1,2]$, $E[\Lambda_0] = \infty$, and for $\tau\in(2,\infty)$
 $E[\Lambda_0] < \infty$ hold by Breiman's theorem; see \cite{JesMik2006}. It is assumed that $P(\Lambda_0 \geq 1) = 1$. 

$L_N = \sum_{i=0}^{N} \Lambda_i$ is the total capacity of the nodes in $G_N$. Using $\Lambda$, we construct a sequence of random graphs $G_N$, $N = 0, 1, \ldots$, where the set of nodes in $G_N$ is $V_N=\{0, \ldots, N\}$. Let $E_N(i, j)$ be the number of edges between the nodes $i$ and $j$ in $G_N$. The graph sequence evolves through the following iterative procedure starting from the initial graph $G_0$ that contains a single isolated node.

\textbf{Step 1.} One forms a random number of self-loops around the initial node in $G_0$ such that $E_0(0,0)$ is a r.v. with   distribution Poisson($\Lambda_0$).

\textbf{Step N + 1.} To construct $G_{N+1}$ from $G_N$, perform the following changes.

(i) New edges are added to the new node $N+1$, namely, for $i \in \{0, \ldots, N + 1\}$, the number $E_{N+1}(i, N+1)$ is generated as a Poisson r.v. with mean
\begin{equation}
\label{eq:e1}
\mathbb{E}[E_{N+1}(i, N + 1)] = \frac{\Lambda_i \Lambda_{N+1}}{L_{N+1}}.
\end{equation}
New edges are added independently of the existing graph structure.

(ii) 
Each old edge of $G_N$ is
independently deleted with probability
\[P_{N+1} := 1 - \frac{L_N}{L_{N+1}}.\]
\par
The evolution model can play a double role:
it serves for the growing network (evolution) and for the spread of information. We focus on spreading a single message among the nodes during the evolution by the Norros-Reittu PAM.
\par At each evolution step $k$, $k=1,2,...$ we may observe two events: either $\#S_k=\#S_{k-1}$ or  $\#S_k=\#S_{k-1}+1$. The increase in  $\#S_{k}$ is considered as success. It is assumed that the initial node has a single message, i.e. $\#S_0=1$.

 Our achievements will be based on the following lemma proved for the linear PAM proposed in \cite{WanWangDavRes}.
\begin{Lemma}\label{lemma1}\cite{Mark2024}
The conditional pmf of  $\#S_{k}$  for a maximum number of evolution steps $K^*$ in fixed time $T^*$ is the following.
For $1\le i\le K^*$
it holds 
\begin{eqnarray*}\label{5b}
P\{\#S_{K^*}=i|G_{K^*-1}\}
&=& e^{-\lambda T^*}\sum_{k=i}^{\infty}\frac{(\lambda T^*)^k}{k!}P\{\#S_k=i|G_{k-1} \},
\end{eqnarray*}
where 
\begin{eqnarray}\label{12a}
&& P\{\#S_k=i|G_{k-1} \}= \sum_{c,k,i-1}\prod_{n=1}^{i-1}p_{j_n}\prod_{m=i}^{k}(1-p_{j_m})\Ii\{k\ge i\ge 2\}\nonumber
\\
&+&
\prod_{m=1}^{k}(1-p_{m})\Ii\{k\ge i=1\}
+\prod_{n=1}^{k}p_{n} \Ii\{i=k+1\}, 
= \psi(i,j,k),
\end{eqnarray}
\begin{eqnarray}\label{12aa}
&& P\{\#S_k=i|G_{k-1} \}=0,\qquad otherwise,
\end{eqnarray}
and  success probabilities $p_1,..., p_k$ correspond to the PAM proposed in \cite{WanWangDavRes}.
\end{Lemma}

 Here, $\sum_{c,k,j}$ denotes the sum of all $\binom{k}{j}=k!/(j!(k-j)!)$ distinct index combinations among $\{j_1, j_2,..., j_k\}$ of length $j$.
 The sum $\sum_{c,k,j}$ corresponds to a Poisson binomial distribution.
 
An integer-valued r.v. $X$ is called the Poisson binomial and is denoted as $X\sim PB(p_{1},p_{2},...,p_{k})$, if $X=^d \xi_1+...+\xi_k$, where $\xi_1,...,\xi_k$ are independent Bernoulli rvs with parameters $p_{1},p_{2},\dots ,p_{k}$. The probability distribution of $X$ is $P\{X=j\}=\sum_{A\in [k],\|A\|=j}\left(\prod_{i\in A}p_i\prod_{i\neq A}(1-p_i)\right)$, where the sum ranges over all subsets of $[k]=\{1,...,k\}$ of size $j$ \cite{Tang}.
 
 Our results are based on Lemma \ref{lemma1}, where $p_k$, $k=1,2,...$ correspond to the Norros-Reittu PAM.

\section{Main results}\label{Sec3}
Let the initial graph $G_0$ contain a single node with a message, i.e. $\#V_0=\#S_0=1$, $\#E_0=0$.  By the Norros-Reittu model,
each new node connects all old nodes existing at the previous evolution step by a Poisson-distributed random number of new edges. The message cannot be transmitted to a new node, if there is no one new edge connecting the latter node with  old nodes with the message. 
In other words, if the number of new edges from a new node to all old nodes with the message is equal to zero, then transmission is unsuccessful.
Otherwise,  the transmission is successful irrespective of the number of new edges. 
\begin{Remark}
The removal of old edges and self-loops does not affect the propagation of the message.
\end{Remark}
In the Norros-Reittu PAM the number of new edges and their appending to all nodes are independent at each evolution step. The distribution (\ref{0}) 
affects the graph structure.

\begin{figure}
    \centering
    \includegraphics[width=1\linewidth]
    {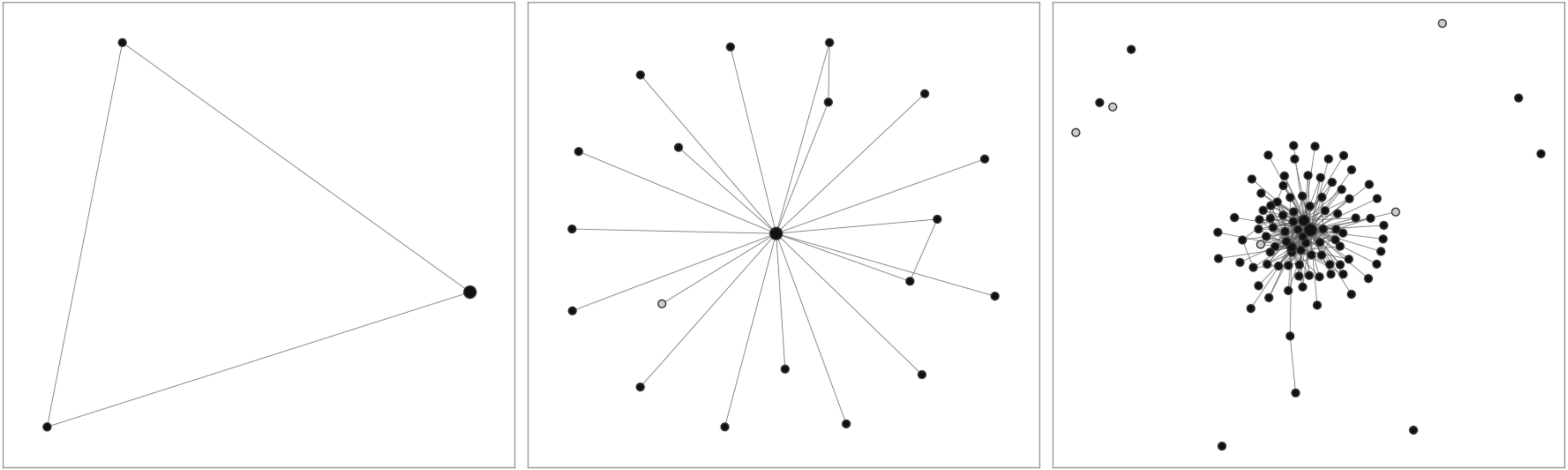}
    \includegraphics[width=1\linewidth]
    {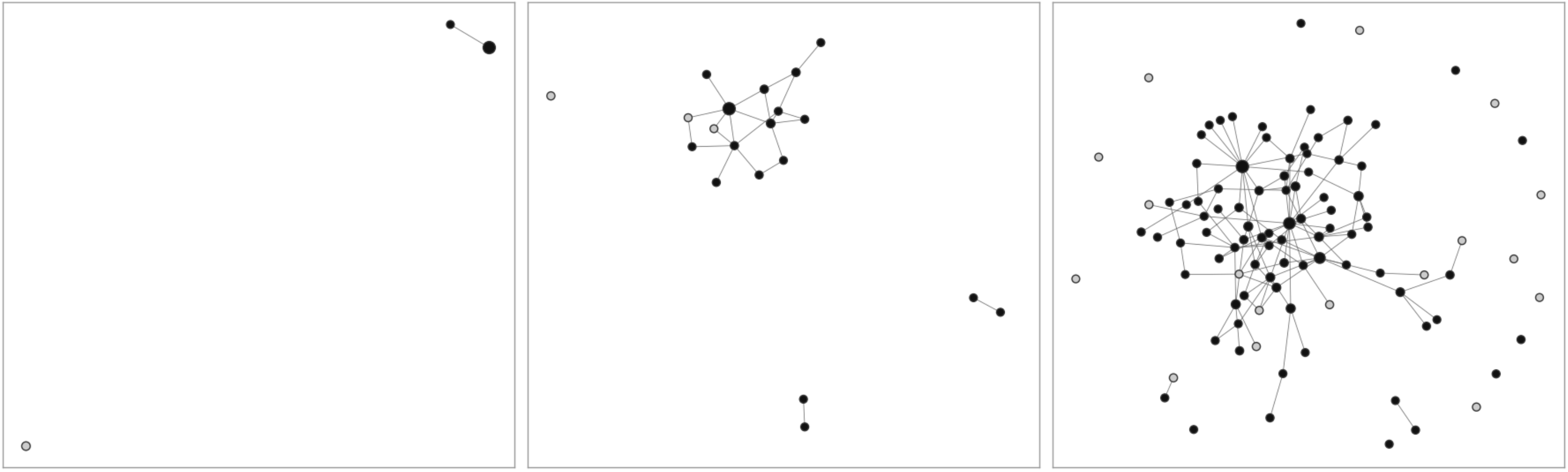}
    \includegraphics[width=1\linewidth]
    {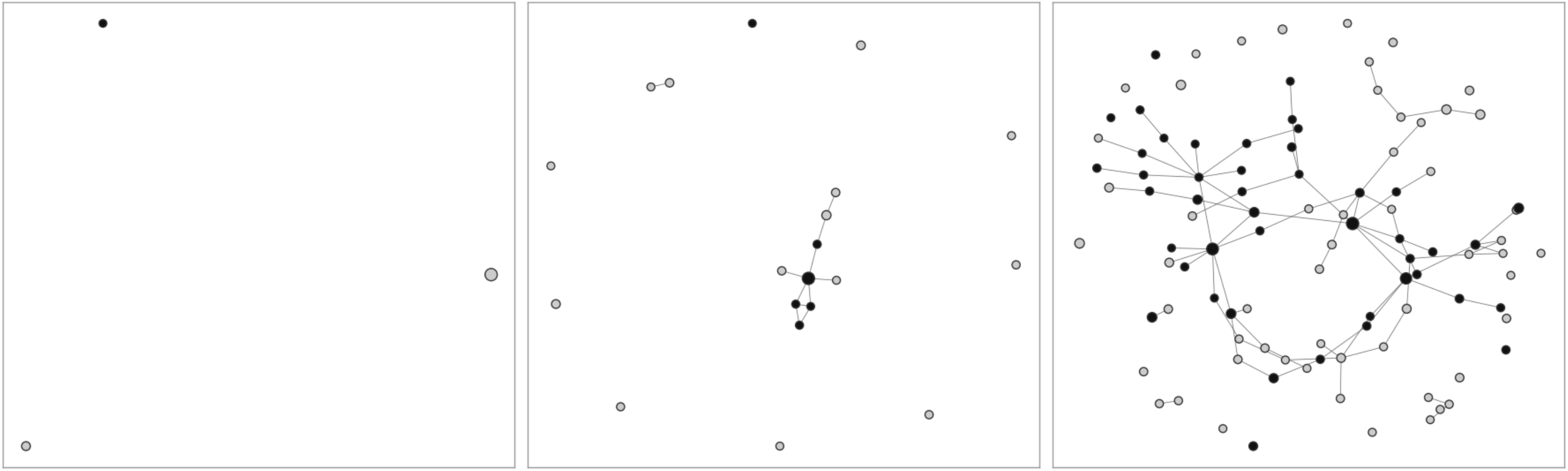}
    \caption{The evolution of a network generated using the Norros-Reittu model with regularly varying distributed node capacities (see (\ref{0}) with $\tau$ = 2.5), where the nodes with messages are marked in black, without message $-$ in grey. Size of circles is proportional to node weight ($\Lambda$). Multiple edges are shown as a single edge, self-loops are not shown. On the left side $-$ the  graph with $N_2=3$; in the middle - Emergent Phase ($N_{19}=20$), where the preferential attachment creates hub structures; on the right $-$ Mature Network ($N_{99}=100$) with the giant connected component. The first line is responsible for $\tau = 1.5$, the second line for $\tau = 2.5$, and the third line for $\tau = 3.5$.}
    \label{fig:graph_evolution}
\end{figure}

As depicted in Figure \ref{fig:graph_evolution}, the graph evolves through distinct topological phases. Early-stage growth ($N \leq 20$) exhibits random tree-like structures, while mature networks ($N \geq 100$) develop scale-free properties consistent with empirical web graphs.

The choice of the distribution parameter \(\tau\) balances heavy-tailed node capacities: the lower \(\tau\) (e.g., \(\tau=1.5\)) leads to numerous high-capacity hubs, enabling rapid message propagation, while the higher \(\tau\) (e.g., \(\tau=3.5\)) results inhomogeneous capacities, slowing the spread.

\subsection{Success probability $p_k$}
\begin{Lemma}\label{lemma2} Under the conditions of graph evolution according to the Norros-Reittu PAM, let 
$\#S_0=1$, and
$\Lambda$ be distributed as in (\ref{0}).
Then the 
success probability $p_{k+1}$, $k=0,1,...$, conditionally on $\Lambda$ and $S_k$ is the following
\begin{equation}
\label{eq:e2}
p_{k+1} = \mathbb{E}\{\#S_{k + 1} - \#S_{k} \mid S_{k}, \Lambda\}
= 1 - \prod_{w \in S_k} e^{-\frac{\Lambda_w \Lambda_{k+1}}{L_{k+1}}}
\end{equation}
\end{Lemma}
\begin{Proof}
By the Norros-Reittu model  it follows that
\begin{itemize}[leftmargin=3.5em, label=$\circ$]
    \item for each node $w \in V_{k}$ number of edges between the new node $k+1$ and an existing node $w$ is defined independently with an intensity determined by their capacities $\Lambda_w$ and $\Lambda_{k+1}$;
    \item the average number of edges between the nodes
    $w\in V_k$ and $k+1$ is given by  \eqref{eq:e1}:
    $$
   \mathbb{E}[E_{k+1}(w, k + 1) \mid \Lambda] = \frac{\Lambda_w \Lambda_{k+1}}{L_{k+1}}, \quad \text{where } L_{k+1} = \sum_{i=0}^{k+1} \Lambda_i.
    $$
\end{itemize}

For the Poisson distribution (\ref{3}), the probability that there are
no edges between nodes $w$ and $k+1$ at step $k+1$ is
$$
P\{E_{k+1}(w, k + 1) = 0 \mid \Lambda\}
= 
e^{-\frac{\Lambda_w \Lambda_{k+1}}{L_{k+1}}}.
$$

Since edges between the new node $k+1$ and 
different existing nodes $w \in S_{k}$ are created independently, the probability that the new node does not connect to any node in $S_{k}$ is equal to the product of probabilities of the absence of edges for all $w \in S_{k}$
$$
P\{S_{k} = S_{k + 1} \mid S_{k}, \Lambda\}
= \prod_{w \in S_{k}} e^{-\frac{\Lambda_w \Lambda_{k+1}}{L_{k+1}}}.
$$
Hence, (\ref{eq:e2}) follows.
\end{Proof}
\begin{Remark} Lemma \ref{lemma2} is valid for any initial graph with at least one node with the message, $\#S_0\ge 1$.
\end{Remark}

\subsection{Probability mass function of $S_k$}
To analyze the pmfs of \(N_k\) and \(\#S_k\), we will use the Poisson binomial distribution model. Each evolution step \(k\) corresponds to an independent experiment defined as follows:  

\begin{enumerate}[leftmargin=3.5em, label=$\circ$]
    \item \textbf{Experiment Description}: At step \(k+1\), a new node \(v_{k+1}\) is added to the graph. Edges between \(v_{k+1}\) and existing nodes are created independently, with the number of edges between \(v_{k+1}\) and each existing node \(w\) following a Poisson distribution with mean \(\frac{\Lambda_w \Lambda_{k+1}}{L_{k+1}}\). This process is independent of previous steps.  

    \item \textbf{Message Transmission Condition}: The message is transmitted to \(v_{k+1}\) \textit{if and only if} at least one edge is created between \(v_{k+1}\) and any node \(w \in S_k\). If no such edges are formed, the transmission fails.  
\end{enumerate}

The success probability \(p_{k+1}\) at step \(k+1\) is given by \eqref{eq:e2}. Since edge creation is independent across steps, the sequence of successes \(\{p_k\}\) forms a Poisson binomial distribution model.





\par
Further analysis is based on the arguments of \cite{Mark2024}.


From Lemma \ref{lemma1} we have that the pmf of $S_k$, $k \geq 1$, for the maximum number of evolution steps  $K^{*}$ within a fixed time $T^{*}$ is defined as follows:
\begin{equation}
\label{eq:11}
P\{\#S_{K^{*}}=i \mid \Lambda\} =e^{-\lambda T^{*}} \sum_{k=i}^{\infty} \frac{(\lambda T^{*})^{k}}{k!} P\{\#S_{k}=i \mid \Lambda\}
\end{equation}
where $P(\#S_k=i \mid \Lambda)$ substitutes
$P\{\#S_k=i|G_{k-1}\}$ in (\ref{12a}), (\ref{12aa}).
\begin{Corollary}\label{corol1}
Let the conditions of Lemma \ref{lemma2} be fulfilled.
The probability of complete non-propagation at evolution step $K \geq1$ is given by
\begin{equation}\label{eq:corol1}
P\{\#S_{K}=1 \mid \Lambda\}= e^{-\alpha}, \qquad \alpha=\Lambda_0 \sum_{k=1}^K \frac{\Lambda_{k}}{L_{k}}.
\end{equation}
\end{Corollary}
\begin{Proof}
As $\#S_k=1$ for each $1 \leq k\leq K$, then by (\ref{eq:e2})
\begin{equation*}
1 - p_k = e^{-{\Lambda_0 \Lambda_{k}} / {L_{k}}}.
\end{equation*}
Hence, (\ref{eq:corol1}) follows by (\ref{12a}). 
\end{Proof}
\par
Now we obtain (\ref{2}) and (\ref{3a}).


\begin{Lemma}\label{lemma4}
Under the conditions of Lemma \ref{lemma2}, 
the pmf of $N_k$, $k \geq 1$, for the maximum number of evolution steps  $K^{*}$ within a fixed time $T^{*}$ is defined as follows:
\begin{eqnarray*}&&
P\left\{N_{K^{*}} = i \right\} = \frac{(\lambda T^{*})^{i-1} e^{-\lambda T^{*}}}{(i-1)!}, \qquad i \geq 1.
\end{eqnarray*}
\end{Lemma}
\begin{Proof}
Consider the process of adding new nodes to a graph.  
Since exactly one node is added at each step, the number of nodes $N_{K^*}$ at step $K^*$ is equal to $K^* + N_0=K^* + 1$.  

According to the condition, the addition of nodes follows a Poisson process, then from the above reasoning and (\ref{3}), the desired result follows
\begin{equation*}
P\left\{N_{K^{*}} = i \right\}
= P\left\{K^{*} = i - 1 \right\}
= \frac{(\lambda T^{*})^{i-1} e^{-\lambda T^{*}}}{(i-1)!}.
\end{equation*}
\end{Proof}
\begin{Lemma}\label{lemma5}
Under the conditions of Lemma \ref{lemma2}, let $K^*$ be determined as in (\ref{3}). Then it holds
\begin{eqnarray*}&&
P\{\#S_{K^*}/N_{K^*} \leq x \mid \Lambda\}
= e^{-\lambda T^{*}} + e^{-\lambda T^{*}} \sum_{k=1}^{\infty}\frac{(\lambda T^{*})^{k}}{k!}\sum_{i=1}^{\lfloor x(k + 1) \rfloor} \psi(i,j,k),
\end{eqnarray*}
where  $\psi(i, j, k)$ is determined by (\ref{12a}).
\end{Lemma}
\begin{Proof}
By (\ref{eq:11}) we have
\begin{equation*}
\begin{aligned}
\label{eq:e5}
&P\{\#S_{K^*}/N_{K^*} \leq x \mid \Lambda\} 
= P\{\#S_{K^*} \leq x(K^* + 1) \mid \Lambda\} \\
&= \sum_{i=1}^{\lfloor x(K^* + 1) \rfloor} P\{\#S_{K^*} = i\mid \Lambda\} \\
&= P\{K^* = 0\} + \sum_{k=1}^{\infty}\sum_{i=1}^{\lfloor x(k + 1) \rfloor} P\{\#S_{K^*} = i, K^* = k \mid \Lambda\} \\
&= e^{-\lambda T^{*}} + e^{-\lambda T^{*}} \sum_{k=1}^{\infty}\sum_{i=1}^{\lfloor x(k + 1) \rfloor}\frac{(\lambda T^{*})^{k}}{k!} P\{\#S_{k}=i \mid \Lambda\}
\end{aligned}
\end{equation*}
Finally, using that $P\{\#S_{k}=i \mid \Lambda\} = \psi(i, j, k)$ (by Lemma \ref{lemma1}) we get
\begin{equation*}
\begin{aligned}
\label{eq:e5}
P\{\#S_{K^*}/N_{K^*} \leq x \mid \Lambda\}
= e^{-\lambda T^{*}} + e^{-\lambda T^{*}} \sum_{k=1}^{\infty}\frac{(\lambda T^{*})^{k}}{k!}\sum_{i=1}^{\lfloor x(k + 1) \rfloor} \psi(i,j,k)
\end{aligned}
\end{equation*}
\end{Proof}
\section{Simulation study}\label{Sec4}
\begin{figure}
    \centering
    \includegraphics[width=1\linewidth]{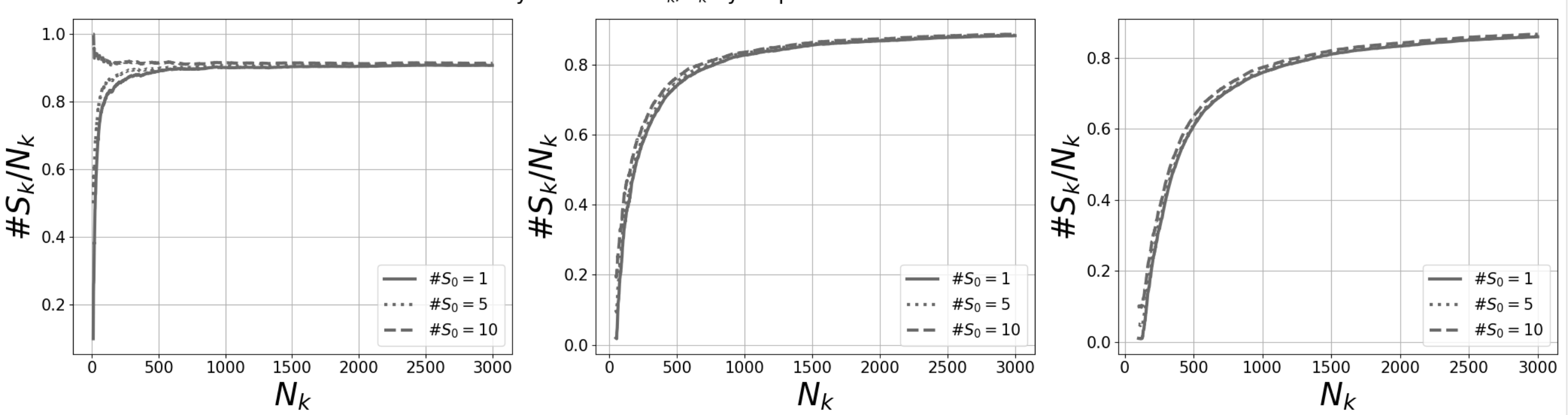}
    \includegraphics[width=1\linewidth]{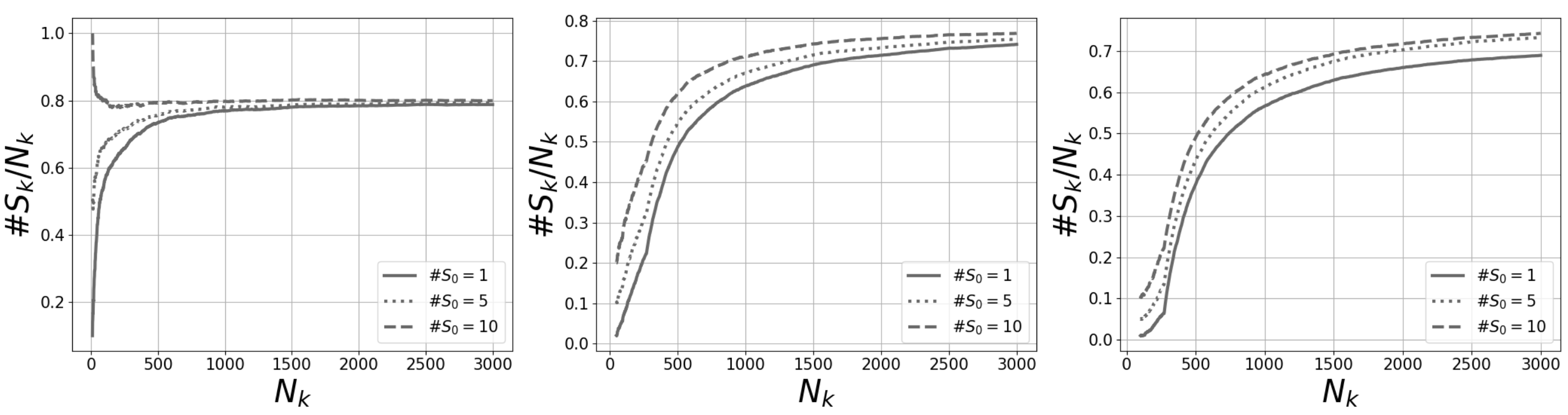}
    \includegraphics[width=1\linewidth]{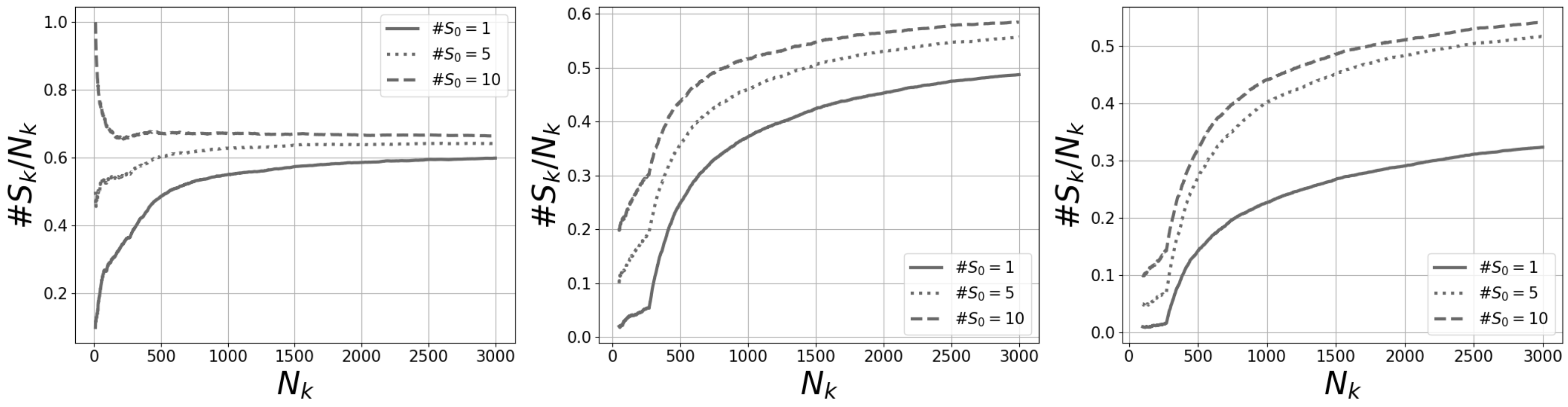}
    \caption{The average value of $\#S_k / N_k$ against $N_k$ for $20$ repetitions of graph evolution by the Norros-Reittu PAM: on the left side $-$ $N_0 = 10$ and $\#S_0 \in \{1, 5, 10\}$; in the middle $-$ $N_0 = 50$ and $\#S_0 \in \{1, 5, 10\}$; on the right $-$ $N_0 = 100$ and $\#S_0 \in \{1, 5, 10\}$, where the first line is responsible for $\tau = 1.5$, the second line for $\tau = 2.5$, and the third line for $\tau = 3.5$.}
    \label{fig:model_norros}
\end{figure}

The simulation results of the Norros-Reittu model were used to plot the dynamics of  
the ratio  
\( \#S_k/N_k \) 
as a function of the size of the  graph \( N_k \) for the following parameters: \textbf{the node capacity distribution} corresponds to  
parameters \( \tau \in \{1.5, 2.5, 3.5 \}\) in (\ref{0}), governing the heavy-tailed capacity distribution; \textbf{initial conditions} corresponding to the initial graph size \( N_0 \in \{10, 50, 100\} \) and the initial number of nodes with the message \( \#S_0 \in \{1, 5, 10\} \), 
see Figure \ref{fig:model_norros}.

Each plot was averaged over $20$ simulation runs. For each \( \tau \), three panels are shown, corresponding to different \(N_0\) with lines representing \( \#S_k/N_k \) dynamics for different \( \#S_0 \). We may conclude the following.

\begin{itemize}[leftmargin=2.5em, label=$\circ$]
    \item \textbf{Impact of initial nodes with the message (\( S_0 \)):} Increasing \( \#S_0 \) accelerates the growth of the \( \#S_k/N_k \) ratio during early graph evolution. However, in case $\tau=1.5$ for large \( N_k \) (\( N_k > 800 \)), the dependence on \( \#S_0 \) diminishes, and all curves converge to a common limit.  
    \item \textbf{Role of parameter \( \tau \):} The lower \( \tau \) values corresponding to distributions (\ref{0}) with heavier tails 
    promote faster message spreading due to high-capacity nodes actively forming connections. For instance, at \( \tau = 1.5 \) and $N_0 = 10$, the \( \#S_k/N_k \) ratio reaches \( 0.7 \)–\( 0.8 \) by \( N_k = 200 \), and ultimately exceeds \( 0.9 \) at \( N_k \approx 800 \). In contrast, for \( \tau = 3.5 \), similar values (\( \#S_k/N_k > 0.9 \)) are not observed even at \( N_k = 3000 \).  
    \item \textbf{Impact of initial graph size (\( N_0 \)):} Larger \( N_0 \) values  amplify the effect of initial conditions: differences between \( \#S_0 = 1, 5, 10 \) are significant. For smaller \( N_0 \) the dynamics rapidly stabilize, and variations in \( S_0 \) become less pronounced.  
    \item \textbf{Process stabilization:} In all cases, the \( \#S_k/N_k \) ratio stabilizes as the graph grows, indicating equilibrium in the message dissemination. 
\end{itemize}

\section{Conclusions}\label{Sec5}
We study the propagation of a message in the Norros-Reittu PAM within a fixed time interval \( T^{*} \). The evolution starts from an initial graph containing the \( S_0 \) nodes with the message, where the network grows according to a heavy-tailed capacity distribution governed by the tail index \( \tau -1\). The message is propogated to a new node only if it connects to at least one existing node with the message, with probabilities determined by the node capacities \( \Lambda_i \).  

We obtained the following results.
\begin{itemize}[leftmargin=2.5em, label=$\circ$]
    \item 
    The pmfs for \( \#S_{K^{*}} \) (the number of nodes with the message), \( N_{K^{*}} \) (the total number of nodes), and the distribution function of their ratio \( \#S_{K^{*}}/N_{K^{*}} \) at time \( T^{*} \).  
    \item Analytical expression for the success probability \( p_k \), linking it to the Poisson binomial distribution and the node capacity dynamics.  
    \item Empirical validation through simulations, demonstrating the critical impact of \( \tau \), \( N_0 \), and \( \#S_0 \) on the \( \#S_k/N_k \) ratio. Evidence of equilibrium in the message dissemination
    as \( N_k \) approaches $3000$ nodes.  
\end{itemize}

\noindent Our future work may concern the following problems:  
\begin{itemize}[leftmargin=2.5em, label=$\circ$]
    \item a generalization to time-dependent parameters (e.g., dynamic of \( \tau \)), node removal mechanisms and different initial graphs; 
    \item a robust estimation of parameters of the PAM associated with the network;
    \item application to real-world networks, such as social media or IoT systems.  
\end{itemize}


\begin{thebibliography}{99}
\bibitem{Norros} Norros, I., Reittu, H.:  On a conditionally poissonian graph process. Advances in Applied Probability, 38(1), 59--75. DOI: 10.1239/aap/1143936140. (2006)

\bibitem{WanWangDavRes} Wan, P.,  Wang, T.,  Davis, R.A.,  Resnick, S.I.:  Are extreme value estimation methods useful for network data?  Extremes,  vol. 23, pp. 171--195. (2020)

\bibitem{Wang} Wang, T., Resnick, S.: Poisson Edge Growth and Preferential Attachment Networks. Methodol Comput Appl Probab 25, 8 (2023).{https://doi.org/10.1007/s11009-023-09997-y}

\bibitem{Mosk-Aoyama2006} Mosk-Aoyama, D., Shah, D.: Computing separable functions via gossip. In Proceedings of the 25th ACM symposium on Principles of distributed computing (PODC '06), ACM, New York, USA  113--122. DOI: 10.1145/1146381.1146401. (2006)

\bibitem{Censor-Hiller2010} Censor-Hillel, K., Shachnai, H.: Partial Information Spreading with Application to Distributed Maximum Coverage. In Proceedings of the 29th ACM symposium on Principles of distributed computing (PODC'10), ACM, New York, USA.  161--170. DOI: 10.1145/1835698.1835739 (2010)

\bibitem{Pat} Patwardhan, S., Rao, V.K., Fortunato, S., Radicchi, F.: Epidemic Spreading in Group-Structured Populations. Physical Review X, 13(4),  041054. DOI: 10.1103/PhysRevX.13.041054 (2023)

\bibitem{Lang} Zeng, L., Tang M., Liu, Y., Yeop Yang, S., Do, Y.: Contagion dynamics in time-varying metapopulation networks with node’s activity and attractiveness.  arXiv:2311.13856v1 [physics.soc-ph] (2023)


\bibitem{Pang} Pang, G.,  Pardoux, E., Velleret, A.: Stochastic SIR model with individual heterogeneity and infection-age dependent infectivity on large non-homogeneous random graphs. arXiv:2502.04225v1 (2025)

\bibitem{Shah} Shah, D.: Gossip algorithms. Foundations and Trends in Networking 3(1),  1--125 (2008)

\bibitem{JesMik2006} Jessen, A.H.,  Mikosch, T.: Regularly varying functions. Publ. Inst. Math. (Beograd) (N.S.), 80, 171--192. (2006)

\bibitem{Mark2024} Markovich, N.M.: Spreading of a limited lifetime information in   networks evolved by preferential attachment.
Submission to e-Journal Reliability: Theory and Applications (2024)

\bibitem{Tang} Tang, W., Tang, F.: The Poisson Binomial Distribution— Old \& New. 
Statist. Sci. 38(1) 108--119. DOI: 10.1214/22-STS852. (2023)

\end{thebibliography}
\end{document}